\documentclass{article}
\usepackage{amsmath}
\usepackage{amsfonts}
\usepackage{amssymb}
\usepackage{makeidx}
\usepackage{graphicx}
\usepackage{natbib}
\usepackage{dsfont}
\usepackage{mathrsfs}

\usepackage{amsthm}

\usepackage{listings,color,relsize}
\newtheorem{theorem}{Theorem}[section]

\theoremstyle{definition}

\newtheorem{example}[theorem]{Example}

\newtheorem{lemma}[theorem]{Lemma}
\newtheorem{prop}[theorem]{Proposition}

\renewcommand{\d}{\,\mathrm{d}}% differential
\renewcommand{\L}{\mathscr{L}} % infinitesimal generator
\newcommand{\img}[1]{\medskip\centering \includegraphics[width=8cm,height=6cm]{eps/#1.eps}}
  \newtheorem{corollary}{Corollary}[section]
  \newtheorem{ex}{Example}[section]
  \newtheorem{rem}{Remark}[section]
  \renewcommand\bibname{References}

\def\bth{\medskip\begin{theorem}}
\def\eth{\end{theorem}\medskip}
\def\brm{\medskip\begin{rem}}
\def\erm{\end{rem}\medskip}
\def\blm{\medskip\begin{lemma}}
\def\elm{\end{lemma}\medskip}
\newcommand{\be}{\begin{equation}}
\newcommand{\ee}{\end{equation}}

\def\bdf{\medskip\begin{defi}}
\def\edf{\end{defi}\medskip}
\def\bco{\medskip\begin{cor}}
\def\eco{\end{cor}\medskip}
\def\bpr{\medskip\begin{prop}}
\def\epr{\end{prop}\medskip}
\def\bcj{\medskip\begin{conj}}
\def\ecj{\end{conj}\medskip}
\def\bex{\medskip\begin{exa}}
\def\eex{\end{exa}\medskip}
\def\bpb{\medskip\begin{prob}}
\def\epb{\end{prob}\medskip}

\begin{document}
%\title{Numerical Algorithm of Backward Doubly Stochastic Differential Equations and SPDEs}

\author{Yufeng Shi \thanks{
Partially supported by National Natural Science Foundation of China Grant
10771122, Natural Science Foundation of Shandong Province of China Grant
Y2006A08 and National Basic Research Program of China (973 Program, No.
2007CB814900). e-mail: yfshi@sdu.edu.cn},
Weiqiang Yang\thanks{Corresponding author, ywq@sdu.edu.cn, http://finance.math.sdu.edu.cn/faculty/yang},
Jing Yuan \\
%EndAName
School of Mathematics, Shandong University \\
Jinan 250100, China}
\title{Numerical Computation for Backward Doubly SDEs and SPDEs}
\maketitle

\begin{abstract}
In this paper we present two numerical schemes of approximating solutions of
backward doubly stochastic differential equations (BDSDEs for short). We
give a method to discretize a BDSDE. And we also give the proof of the
convergence of these two kinds of solutions for BDSDEs respectively. We give
a sample of computation of BDSDEs.
\end{abstract}



\maketitle
%\frontmatter   % for book
%\tableofcontents

%\mainmatter

\section{Introduction}
\label{chap:chap1}
%%%%%%%%%%%%%%%%%%%%%%%%%%%%%%%%%%%%%%%%%%%%%%%%%%%%%%%
%%%%%%%%%%%%%%%%%%%%%%%%%%%%%%%%%%%%%%%%%%%%%%%%%%%%%%%
%%%%%%%%%%%%%%%%%%%%%%%%%%%%%%%%%%%%%%%%%%%%%%%%%%%%%%%
Since Pardoux and Peng introduced backward stochastic differential
equation (BSDE), the theory of which has been widely used and
developed, mainly because of a large part of problems in
mathematical finance can be treated as a BSDE. However it is known
that only a limited number of BSDE can be solved explicitly. To
develop numerical method and numerical algorithm is very helpful,
theoretically and practically. Recently many different types of
discretization of BSDE and the related numerical analysis were
introduced.

On the other hand, Paroux and Peng \cite{PP} introduced a new class
of backward stochastic differential equations-backward“doubly”
stochastic differential equations and also showed the existence and
uniqueness of the solution of BDSDE. But until now little work is
devoted to the numerical method and the related numerical analysis.
Here following the approach of M\'{e}min, Peng and Xu \cite{MPX}, we
present two numerical schemes of approximating solutions of BDSDE,
and proved the convergence of these two kinds of solutions for
BDSDEs, respectively. First of the proofs makes use of and extends
Donsker-Type theorem.

This paper is organized as follows. In section 2, we introduce some
fundamental knowledge and assumptions of BDSDEs. In section 3, the
discrete BDSDE and solutions are presented. In section 4, we will
give our main results: the proof of convergence of numerical
solutions for BDSDEs in two different schemes.

\section{Some Preliminaries}
\label{chap:chap2}

%%%%%%%%%%%%%%%%%%%%%%%%%%%%%%%%%%%%%%%%%%%%%%%%%
%%%%%%%%%%%%%%%%%%%%%%%%%%%%%%%%%%%%%%%%%%%%%%%%%
%%%%%%%%%%%%%%%%%%%%%%%%%%%%%%%%%%%%%%%%%%%%%%%%%

Let ($\Omega, \mathcal{F}, \mathbf{P}$) be a complete probability
space, and $T>0$ be fixed throughout this paper. Let
\{$\mathit{W}_{t}$, $0\leq t\leq T$\}
 and \{$\mathit{B}_{t}$, $0\leq t\leq T$\} be two mutually independent
standard Brownian motion processes, with values respectively in
$\mathbb{R}^{d}$and in $\mathbb{R}^{l}$, define on
($\Omega$,$\mathcal{F}$,$ \mathbf{P} $). For each $t\in [0,T]$, we
define $$\mathcal{F}_{t}\doteq \mathcal{F}_{t}^{W}\vee
\mathcal{F}_{t,T}^{B}$$ where for any process \{$\eta _{t}$\},
$\mathcal{F}_{s,t}^{\eta }=\sigma \{\eta _{r}-\eta _{s};s\leq r\leq
t\}$, $\mathcal{F}_{t}^{\eta }=\mathcal{F} _{0, t}^{\eta }$.

For any $n\in \mathbb{N}$, let $\mathcal{M}^{2}(0,T;\mathbb{R}^{n})$
denote the set of (classes of $d\mathbf{P}\times dt$ $a.e.$ equal)
$n$ dimensional jointly measurable random processes \{$\varphi
_{t};t\in [0,T]$\} which satisfy:

(i). $\mathbb{E}\int\limits_{0}^{T}\left\vert \varphi_{t}\right\vert
^{2}dt<\infty $

(ii). $\varphi _{t}$ is $\mathcal{F}_{t}$-measurable, for $a.e.t\in
\lbrack 0,T].$

We denote similarly by $\mathcal{S}^{2}([0,T];\mathbb{R}^{n})$ the
set of continuous $n$ dimensional random processes which satisfy:

(i). $\mathbb{E}(\sup_{0\leq t\leq T}\left\vert \varphi
_{t}\right\vert ^{2})<\infty $

(ii). $\varphi _{t}$ is $\mathcal{F}_{t}$-measurable, for any
$t\in[0,T].$

Let $$f:\Omega \times \lbrack 0,T]\times \mathbb{R}^{k}\times
\mathbb{R}^{k\times d}\rightarrow \mathbb{R}^{k}$$
$$g:\Omega \times \lbrack 0,T]\times \mathbb{R}^{k}\times
\mathbb{R}^{k\times d}\rightarrow \mathbb{R}^{k\times l}$$ be
jointly measurable and such that for any $(y,z)\in
\mathbb{R}^{k}\times \mathbb{R}^{k\times d},$
$$f(\cdot ,y,z)\in \mathcal{M}^{2}(0,T;\mathbb{R}^{k})$$
$$g(\cdot ,y,z)\in \mathcal{M}^{2}(0,T;\mathbb{R}^{k\times l})$$

We assume moreover that there exist constants $K>0$ and $0<\alpha
<1$ such that for any $(\omega ,t)\in \Omega
\times[0,T],(y_{1},z_{1}), (y_{2},z_{2})\in \mathbb{R}^{k}\times
\mathbb{R}^{k\times d},$

$$(H.1)\ \
\left\vert f(t,y_{1},z_{1})-f(t,y_{2},z_{2})\right\vert \leq
K(\left\vert y_{1}-y_{2}\right\vert +\left\Vert
z_{1}-z_{2}\right\Vert )$$
$$\hspace{13mm}\left\Vert
g(t,y_{1},z_{1})-g(t,y_{2},z_{2})\right\Vert \leq K\left\vert
y_{1}-y_{2}\right\vert +\alpha \left\Vert z_{1}-z_{2}\right\Vert$$

Given $\xi \in L^{2}(\Omega
,\mathcal{F}_{T},\mathbf{P};\mathbb{R}^{k})$,  we consider the
following backward doubly stochastic differential equation:
$$Y_{t}=\xi +\int\limits_{t}^{T}f(s,Y_{s},Z_{s})ds+\int
\limits_{t}^{T}g(s,Y_{s},Z_{s})dB_{s}-\int\limits_{t}^{T}Z_{s}dW_{s},
\qquad 0\leq t\leq T.$$

We note that the integral with respect to $\{B_{t}\}$ is a "backward
It\^{o} integral" and the integral with respect to $\{W_{t}\}$ is a
standard forward forward It\^{o} integral, see Nualart and Pardoux
\cite{NP}.

Here we mainly study the case when Brownian motion is
one-dimensional. Now we consider the following 1-dimensional BDSDE

\begin{equation}
\label{equation:BDSDE} Y_{t}=\xi
+\int\limits_{t}^{T}f(s,Y_{s},Z_{s})ds+\int
\limits_{t}^{T}g(s,Y_{s},Z_{s})dB_{s}-\int\limits_{t}^{T}Z_{s}\textrm{d}W_{s}
\end{equation}
and the terminal condition is $y_{T}=\xi =\Phi (W_{T})$, where $\Phi
(\cdot ) $ is a functional of Brownian motion
\{$(B_{s},W_{s})_{0\leq s\leq T}$\}, such that $\xi \in
L^{2}(\mathcal{F}_{T})$. Particularly, if $f(\cdot)$, $g(\cdot)$ are
not relative to $t$, $(\ref{equation:BDSDE})$ changes into:
\begin{equation}
Y_{t}=\xi +\int\limits_{t}^{T}f(Y_{s},Z_{s})ds+\int
\limits_{t}^{T}g(Y_{s},Z_{s})dB_{s}-\int\limits_{t}^{T}Z_{s}\textrm{d}W_{s}
\end{equation}

\section{Numerical Scheme of Standard BDSDE}
\label{chap:chap3}

%%%%%%%%%%%%%%%%%%%%%%%%%%%%%%%%%%%%%%%%%%%%%%%
%%%%%%%%%%%%%%%%%%%%%%%%%%%%%%%%%%%%%%%%%%%%%%%
%%%%%%%%%%%%%%%%%%%%%%%%%%%%%%%%%%%%%%%%%%%%%%%
\subsection{The Structure of Numerical Solution}

When $n\in \mathbb{N}$ is big enough, we divide the time interval
$[0,T]$ into $ n \, $parts: $0=t_{0}<t_{1}<\cdot \cdot \cdot
<t_{n}=T$, $\delta :=t_{j}-t_{j-1}=\frac{T}{n}$, for $1\leq j\leq
n$.

Now we define the scaled random walk $\{B_{\cdot }^{n}\},\{W_{\cdot
}^{n}\}$ , by setting $B_{0}^{n}=W_{0}^{n}=0,$

$$B_{t}^{n}=\sqrt{\delta }\sum_{j=1}^{[t/\delta ]}\varepsilon
_{j}^{n},W_{t}^{n}=\sqrt{\delta }\sum_{j=1}^{[t/\delta ]}\beta
_{j}^{n},\qquad 0\leq t\leq T$$ where $\{\varepsilon
_{j}^{n}\}_{j=1}^{n},\{\beta _{j}^{n}\}_{j=1}^{n}$ are two mutually
independent Bernoulli sequences, which are i.i.d. random variable
satisfying
$$
\varepsilon _{m}^{n}=\beta _{r}^{n}=\left\{\begin{array}{rl}
+1,p=0.5\\ -1,p=0.5\end{array}\right.
$$

Obviously, $B_{t}^{n},$ $W_{t}^{n}$ are both
$\mathcal{F}_{t}$-measurable processes who take discrete values,
denote $ B_{j}^{n}=B_{t_{j}}^{n},W_{j}^{n}=W_{t_{j}}^{n}$, we get
$B_{j}^{n}=\sqrt{\delta }\sum\limits_{m=1}^{j}\varepsilon
_{m}^{n},W_{j}^{n}=\sqrt{\delta } \sum\limits_{r=1}^{j}\beta
_{r}^{n}$. And we define the discrete filtrations $
\mathcal{G}_{j}^{n}=\sigma \{\beta _{1},...,\beta _{j}\}=\sigma
\{W_{t}^{n},0\leq t\leq t_{j}\},\mathcal{G}_{jj}^{n}=\sigma \{\beta
_{1},...,\beta _{j}\}\vee \sigma \{\varepsilon
_{j+1},...,\varepsilon _{n}\}=\sigma \{W_{t}^{n},0\leq t\leq
t_{j}\}\vee \sigma \{B_{t}^{n},t_{j+1}\leq t\leq T\},$

Then, on the small interval $[t_{j},t_{j+1}]$, the equation

\begin{equation}
Y_{t_{j}}=Y_{t_{j+1}}+\int\limits_{t_{j}}^{t_{j+1}}f(s,Y_{s},Z_{s})ds+\int
\limits_{t_{j}}^{t_{j+1}}g(s,Y_{s},Z_{s})dB_{s}-\int
\limits_{t_{j}}^{t_{j+1}}Z_{s}dW_{s}
\end{equation}
can be approximated by the discrete equation
\begin{equation}
y_{j}^{n}=y_{j+1}^{n}+f(t_{j},y_{j}^{n},z_{j}^{n})\delta
+g(t_{j+1},y_{j+1}^{n},z_{j+1}^{n})(B_{j+1}^{n}-B_{j}^{n})-z_{j}^{n}(W_{j+1}^{n}-W_{j}^{n})
\end{equation}
i.e.
\begin{equation}
y_{j}^{n}=y_{j+1}^{n}+f(t_{j},y_{j}^{n},z_{j}^{n})\delta
+g(t_{j+1},y_{j+1}^{n},z_{j+1}^{n})\sqrt{\delta }\varepsilon
_{j+1}-z_{j}^{n} \sqrt{\delta }\beta _{j+1}
\end{equation}

For sake of simplicity, here we just consider the situation in which
$f,g$ are not relative to $t$.

\blm
Let $y_{j+1}^{n}$ be a given
$\mathcal{G}_{j+1,j+1}^{n}$-measurable random variable. Then, when$\
\delta <1/k$, there exists a unique $
\mathcal{G}_{jj}^{n}$-measurable pair ($y_{j}^{n},z_{j}^{n}$)
satisfying the equation:
\begin{equation}\label{equation:discreteBDSDE}
y_{j}^{n}=y_{j+1}^{n}+f(y_{j}^{n},z_{j}^{n})\delta
+g(y_{j+1}^{n},z_{j+1}^{n})\sqrt{\delta }\varepsilon
_{j+1}-z_{j}^{n}\sqrt{ \delta }\beta _{j+1}
\end{equation}
\elm
\textbf{Proof}. We set $Y_{j+1}^{+}=y_{j+1}^{n}\mid _{\beta
_{j+1}=1},Y_{j+1}^{-}=y_{j+1}^{n}\mid _{\beta
_{j+1}=-1},y_{j}^{+}=y_{j}^{n}\mid _{\varepsilon
_{j+1}=1},y_{j}^{-}=y_{j}^{n}\mid _{\varepsilon _{j+1}=-1}.$

Both $Y_{j+1}^{+},Y_{j+1}^{-}$ are
$\mathcal{G}_{j+1,j+1}^{n}$-measurable. Then equation
$(\ref{equation:discreteBDSDE})$ is equivalent to the following
algebraic equations:
$$y_{j}^{+}=Y_{j+1}^{+}+f(y_{j}^{+},z_{j}^{+})\delta
+g(Y_{j+1}^{+},z_{j+1}^{+})\sqrt{ \delta }-z_{j}^{+}\sqrt{\delta }$$

$$y_{j}^{+}=Y_{j+1}^{-}+f(y_{j}^{+},z_{j}^{+})\delta
+g(Y_{j+1}^{-},z_{j+1}^{-})\sqrt{ \delta }+z_{j}^{+}\sqrt{\delta }$$

$$y_{j}^{-}=Y_{j+1}^{+}+f(y_{j}^{-},z_{j}^{-})\delta
-g(Y_{j+1}^{+},z_{j+1}^{+})\sqrt{ \delta }-z_{j}^{-}\sqrt{\delta }$$

$$y_{j}^{-}=Y_{j+1}^{-}+f(y_{j}^{-},z_{j}^{-})\delta
-g(Y_{j+1}^{-},z_{j+1}^{-})\sqrt{ \delta }+z_{j}^{-}\sqrt{\delta }$$

Solving these equations, we can get
$$z_{j}^{+}=\frac{1}{2\sqrt{\delta
}}(Y_{j+1}^{+}-Y_{j+1}^{-})+\frac{1}{2}
[g(Y_{j+1}^{+},z_{j+1}^{+})-g(Y_{j+1}^{-},z_{j+1}^{-})]$$

$$z_{j}^{-}=\frac{1}{2\sqrt{\delta
}}(Y_{j+1}^{+}-Y_{j+1}^{-})-\frac{1}{2}
[g(Y_{j+1}^{+},z_{j+1}^{+})-g(Y_{j+1}^{-},z_{j+1}^{-})]$$

$$y_{j}^{+}-f(y_{j}^{+},z_{j}^{+})\delta
=\frac{1}{2}(Y_{j+1}^{+}+Y_{j+1}^{-})+\frac{\sqrt{ \delta
}}{2}[g(Y_{j+1}^{+},z_{j+1}^{+})+g(Y_{j+1}^{-},z_{j+1}^{-})]$$

$$y_{j}^{-}-f(y_{j}^{-},z_{j}^{-})\delta
=\frac{1}{2}(Y_{j+1}^{+}+Y_{j+1}^{-})-\frac{\sqrt{ \delta
}}{2}[g(Y_{j+1}^{+},z_{j+1}^{+})+g(Y_{j+1}^{-},z_{j+1}^{-})]$$

That is to say:
$$z_{j}^{n}=\frac{1}{2\sqrt{\delta
}}(Y_{j+1}^{+}-Y_{j+1}^{-})+\frac{1}{2}
[g(Y_{j+1}^{+},z_{j+1}^{+})-g(Y_{j+1}^{-},z_{j+1}^{-})]\varepsilon_{j+1}$$

$$y_{j}^{n}-f(y_{j}^{n},z_{j}^{n})\delta
=\frac{1}{2}(Y_{j+1}^{+}+Y_{j+1}^{-})+\frac{\sqrt{ \delta
}}{2}[g(Y_{j+1}^{+},z_{j+1}^{+})+g(Y_{j+1}^{-},z_{j+1}^{-})]\varepsilon
_{j+1}$$

We can simulate a sample path of $\{\varepsilon_{j}\}$, then we
calculate the corresponding BSDE along with the sequence. It is
indeed a kind of Monte-Carlo method.
\begin{ex}
If $f(y, z) = ay + bz$,
$$Y^n_{j} =\frac{\frac{1}{2}(Y_{j+1}^{+}+Y_{j+1}^{-})+\frac{\sqrt{ \delta
}}{2}[g(Y_{j+1}^{+},Z_{j+1}^{+})+g(Y_{j+1}^{-},Z_{j+1}^{-})]\varepsilon
_{j+1}+bZ^n_{j}\delta}{1-a\delta}.$$
\end{ex}

\vspace{10mm} The calculation begins at the terminal time $t_{n}=T$,
with $ y_{n}^{n}=\xi ^{n}$, which is given and the problem is how to
determine $Z_{n}$. Here we choose the way of setting $Z_T = \nabla
Y_T,$  i.e. $Z_T =\partial_{x}Y_T$.

\begin{ex}
For simplicity, we suppose a linear type, $f(y,z)=0$,
$g(y,z)=ay+bz$, $Y_T=W_T$, then $Z_T=\nabla Y_T$.

We have
$$Y_{n-1}^{+}=Y_{n}^{+}+(aY_{n}^{+}+bZ_{n}^{+})\sqrt{ \delta }-Z_{n-1}^{+}\sqrt{\delta }$$

$$Y_{n-1}^{+}=Y_{n}^{+}+(aY_{n}^{+}+bZ_{n}^{+})\sqrt{ \delta }+Z_{n-1}^{+}\sqrt{\delta }$$

$$Y_{n-1}^{-}=Y_{n}^{+}+(aY_{n}^{+}+bZ_{n}^{+})\sqrt{ \delta }-Z_{n-1}^{-}\sqrt{\delta }$$

$$Y_{n-1}^{-}=Y_{n}^{+}+(aY_{n}^{+}+bZ_{n}^{+})\sqrt{ \delta }+Z_{n-1}^{-}\sqrt{\delta }$$
i.e.
$$Y_{n-1}^{+}=\frac{Y_{n}^{+}+Y_{n}^{-}}{2}(1+a\sqrt{ \delta })+\frac{b}{2}(Z_{n}^{+}+Z_{n}^{-})\sqrt{\delta}$$

$$Z_{n-1}^{+}=\frac{Y_{n}^{+}-Y_{n}^{-}}{2\sqrt{\delta
}}(1+a\sqrt{ \delta })+\frac{b}{2}(Z_{n}^{+}-Z_{n}^{-})$$

$$Y_{n-1}^{-}=\frac{Y_{n}^{+}+Y_{n}^{-}}{2}(1-a\sqrt{ \delta })+\frac{b}{2}(Z_{n}^{+}+Z_{n}^{-})\sqrt{\delta}$$

$$Z_{n-1}^{-}=\frac{Y_{n}^{+}+Y_{n}^{-}}{2\sqrt{\delta
}}(1-a\sqrt{ \delta })+\frac{b}{2}(Z_{n}^{+}-Z_{n}^{-})$$
\end{ex}
\vspace{10mm} After $Z_n$ is calculated, $Y_{j}$ and $Z_{j}$ can be
backwardly step by step, following the way mentioned above.

On the other hand, taking conditional expectation on
$(\ref{equation:discreteBDSDE})$, it follows that
$$z_{j}^{n}=\frac{1}{\sqrt{\delta }}\mathbb{E}[y_{j+1}^{n}\beta
_{j+1}\mid
\mathcal{G}_{jj}^{n}]+\mathbb{E}[g(y_{j+1}^{n},z_{j+1}^{n})\beta
_{j+1}\mid \mathcal{G}_{jj}^{n}]\varepsilon_{j+1}$$

$$y_{j}^{n}-f(y_{j}^{n},z_{j}^{n})\delta =\mathbb{E}[y_{j+1}^{n}\mid
\mathcal{G }_{jj}^{n}]+\sqrt{\delta
}\mathbb{E}[g(y_{j+1}^{n},z_{j+1}^{n})\mid \mathcal{
G}_{jj}^{n}]\varepsilon _{j+1}$$

At the terminal time $t_{n}=T$, consider the mapping $\Psi
(z)=z-\frac{1}{2} [g(Y_{+},z)-g(Y_{-},z)]\varepsilon_{j+1}$, from
the property of $g$, we obtain that the derivative of $\Psi (z)$ on
$z$ is $1$, which implies that the mapping $\Psi (z)$ is a monotonic
mapping. So there exists a unique value $z_{n-1}$ s.t.
$z_{n-1}=\frac{1}{2\sqrt{\delta }}
(Y^{+}-Y^{-})+\frac{1}{2}[g(Y^{+},z_{n-1})-g(Y^{-},z_{n-1})]$ holds.
Consider the mapping $\Theta (y)=y-f(y,z_{j}^{n})\delta $ from the
Lipschitz property of $f$, we obtain
$$\left\langle \Theta
(y)-\Theta (y^{\prime }),y-y^{\prime }\right\rangle \geq (1-\delta
K)\left\vert y-y^{\prime }\right\vert ^{2}>0$$ which implies that
the mapping $\Theta (y)$ is a monotonic mapping. So there exists a
unique value $y$ s.t. $\Theta (y)=\mathbb{E}[y_{j+1}^{n}\mid
\mathcal{G}_{jj}^{n}]+\sqrt{\delta
}\mathbb{E}[g(y_{j+1}^{n},z_{j+1}^{n}) \mid
\mathcal{G}_{jj}^{n}]\varepsilon _{j+1}$ holds, i.e.
$y_{j}^{n}=\Theta ^{-1}(\mathbb{E}[y_{j+1}^{n}\mid
\mathcal{G}_{jj}^{n}]+\sqrt{\delta }\mathbb{
E}[g(y_{j+1}^{n},z_{j+1}^{n})\mid \mathcal{G}_{jj}^{n}]\varepsilon
_{j+1})$.

\textbf{Remark a}. The existence of the solution of discrete BDSDE
only depends on the Lipschitz condition of $f$ on $y$. In fact, if
$f$ does not depend on $y$, we can easily get $\Theta
^{-1}(y)=y+f(z_{j}^{n})\delta ,$ And very obviously, if $g$ does not
depend on $z$, $\{z_{\cdot }^{n}\}$ can be also easily got.

\textbf{Remark b.} In general, if $f$ nonlinearly depends on $y$,
then $ \Theta (y)$ can not be solved explicitly, so sometimes we can
introduce the following scheme, we set
$\overline{Y}_{T}^{n}=\overline{y}_{n}^{n}=\xi ^{n}, $ and starting
from $j=n-1,$ backwardly solve
\begin{equation}\label{equation:modifieddiscreteBDSDE}
\overline{y}_{j}^{n}=\overline{y}_{j+1}^{n}+f(\mathbb{E}[\overline{y}
_{j+1}^{n}\mid \mathcal{G}_{jj}^{n}],\overline{z}_{j}^{n})\delta +g(
\overline{y}_{j+1}^{n},\overline{z}_{j+1}^{n})\sqrt{\delta
}\varepsilon _{j+1}-\overline{z}_{j}^{n}\sqrt{\delta }\beta _{j+1}
\end{equation}
or equivalently,
\begin{equation}
\overline{z}_{j}^{n}=\frac{1}{\sqrt{\delta
}}\mathbb{E}[\overline{y}_{j+1}^{n}\beta _{j+1}\mid
\mathcal{G}_{jj}^{n}]+\mathbb{E}[g(\overline{y}_{j+1}^{n},\overline{z}
_{j+1}^{n})\beta _{j+1}\mid \mathcal{G}_{jj}^{n}]\varepsilon_{j+1}
\end{equation}
\begin{equation}
\overline{y}_{j}^{n}=\mathbb{E}[\overline{y}_{j+1}^{n}\mid
\mathcal{G} _{jj}^{n}]+f(\mathbb{E}[\overline{y}_{j+1}^{n}\mid
\mathcal{G}_{jj}^{n}], \overline{z}_{j}^{n})\delta
+\mathbb{E}[g(\overline{y}_{j+1}^{n},\overline{z}
_{j+1}^{n})]\sqrt{\delta }\varepsilon _{j+1}
\end{equation}
to approximate the solution of $\Theta
(y)=\mathbb{E}[y_{j+1}^{n}\mid \mathcal{G}_{jj}^{n}]+\sqrt{\delta
}\mathbb{E}[g(y_{j+1}^{n},z_{j+1}^{n}) \mid
\mathcal{G}_{jj}^{n}]\varepsilon _{j+1}.$

%%%%%%%%%%%%%%%%%%%%%%%%%%%%%%%%%%%%%%%%%%%%%%%  ywq

\subsection{Monte-Carlo Method}
For Forward-Backward SDEs,
\begin{eqnarray}
  \d X_s &=& b(s,X_s)\d s + \sigma(s,X_s) \d W_s. \\
  X_s &=& x, \ \ \ \ 0\leq s \leq t. \nonumber
\end{eqnarray}
\begin{eqnarray}
  -\d Y_s=f(s,X_s^{t,x}, Y_s, Z_s) \d s- Z_s \d W_s,\\
  Y_T = \Psi(X_T^{t,x}). \nonumber
\end{eqnarray}
where $b$, $\sigma$, $f$ and $\Phi$ satisfy usual assumption.
by Theorem 4.1 of \cite{EPQ}
There exist two function $u(t,x)$ and $d(t,x)$,
such that the solution $(Y^{t,x}, Z^{t,x})$ of BSDE is
\[
  Y_s^{t,x}=u(s,X_s^{t,x}), \ \ \ Z_s^{t,x}=\sigma(s, X_s^{t,x}) d(s,X_s^{t,x}),
  \ \ \ \ t\leq s \leq T, \d P \otimes \d s \textrm{ a.s.}
\]
%\begin{lemma}[Theorem 4.1 of EPQ 1997]
%There exist two $\B([0,T])\otimes\B_e(\R)$ measurable deterministic function $u(t,x)$ and $d(t,x)$,
%such that the solution $(Y^{t,x}, Z^{t,x})$ of BSDE is
%\[
%  Y_s^{t,x}=u(s,X_s^{t,x}), \ \ \ Z_s^{t,x}=\sigma(s, X_s^{t,x}) d(s,X_s^{t,x}),
%  \ \ \ \ t\leq s \leq T, \d P \otimes \d s \textrm{ a.s.}
%\]
%\end{lemma}
The solution of the BSDE is said to be Markovian.
So it is naturally to solve the equation based on a binomial tree of $X_s$.
%Markovian BSDEs is a big class of equation,

%\begin{remark}%The lemma doesn't
%Suppose $X_s =\int_0^s r\d W_r$,
%it isn't a function of $(s, W_s)$,
%\end{remark}

\begin{example}
If $X_s \equiv W_s$, the solution of BSDE is
\[
  Y_s=u(s,W_s), \ \ \ Z_s=d(s,W_s),
  \ \ \ \ t\leq s \leq T, \d P \otimes \d s \textrm{ a.s.}
\]
%The solution can be written in the form of $Y_t=\Phi(t,W_t)$, $Z_t=\Psi(t,W_t)$.
\end{example}

%\begin{remark}
As for BDSDE, the structure of BDSDE is different from BSDE,
%Even though $f$ and $g$ are deterministic,
that the solution is not generally in the form of $Y_t=\phi(t,W_t,B_T-B_t)$, $Z_t=\psi(t,W_t,B_T-B_t)$,
even though $f$ and $g$ are deterministic functions.
%\end{remark}
\begin{example}
\begin{eqnarray*}
  -\d Y_t &= &t\d B_t -Z_t \d W_t,\\
  Y_T     &= &0.
\end{eqnarray*}

\end{example}
The solution is
\begin{eqnarray*}
Y_t &=&\int_t^T s\d B_s,\\
Z_t &=&0.
\end{eqnarray*}

Therefore, $Y_t$ is path dependent on $B.$,
%nodes is doubled on every step, that is to say, it will grow exponentially.
So it's impossible to solve the solution on the nodes of the coupled binomial trees.

If we simulate %a Bernoulli sequence $\{\varepsilon_i\}$ to approximation $\Delta B_t$,
a sample path of $B_t$,
it becomes a classical numerical scheme of BSDE follow the path, which is indeed a Monte-Carlo method,
and the solution surface will vibrate with the sample path of $B_t$.
Yang \cite{Y} gives some comparison examples of numerical solutions and explicit solutions.
\footnote{The reason we don't include the examples in this paper is that arXiv reject figures of large size.}

%so we can generate many sample paths of $B$, and calculate the corresponding
%BSDE each time.
%这样可以去掉原来算法的限制。

%\begin{eqnarray*}
%Y_{i} &=&Y_{i+1}^{+}+f(t_{i},Y_{i},Z_{i})\Delta
%+g(t_{i+1},Y^{+}_{i+1},Z^{+}_{i+1})\sqrt{\Delta }\varepsilon_i -Z_{n-1}\sqrt{\Delta } \\
%Y_{i} &=&Y_{i+1}^{-}+f(t_{i},Y_{i},Z_{i})\Delta
%+g(t_{i+1},Y_{i+1}^{-},Z_{i+1}^{-})\sqrt{\Delta }\varepsilon_i +Z_{n-1}\sqrt{\Delta }
%\end{eqnarray*}%

%We can plot the figure of $Y_0$ against $t_0$, $B_T-B_0$,
%but $Y_0$ is not surely a function of $B_T-B_0$.

%二叉树方法和高斯分布的方法有本质差别吗？二叉树方法应该比较稳定，而高斯分布会不稳定？

%A kind of Monte-Carlo method is
%\begin{itemize}
%  \item Monte-Carlo method to generate sample paths of $B_t$.
%  \item Solve the BSDE associated with the generated $B_t$.
%\end{itemize}

\subsection{Associated SPDE}
For each $\left( t,x\right) \in R^{+}\times R^{d},$ let $\left\{
X_{s}^{t,x},t\leq s\leq T\right\} $ be the solution of the SDE:%
\[
X_{s}^{t,x}=x+\int\limits_{t}^{s}b(X_{r}^{t,x})dr+\int\limits_{t}^{s}\sigma
(X_{r}^{t,x})dW_{r},{\ \ \ }t\leq s\leq T.
\]

\bigskip The following BDSDE%

\[
Y_{s}^{t,x}=h\left( X_{T}^{t,x}\right)
+\int\limits_{s}^{T}f(X_{r}^{t,x},Y_{r}^{t,x},Z_{r}^{t,x})dr+\int%
\limits_{s}^{T}g(X_{r}^{t,x},Y_{r}^{t,x},Z_{r}^{t,x})dB_{r}-\int%
\limits_{s}^{T}Z_{r}^{t,x} dW_{r},{\ \ \ }t\leq s\leq T.
\]

Under assumption (H1) has a unique solution $%
(Y_{s}^{t,x},Z_{s}^{t,x})$, and under some suitable conditions,
\[
u(t,x)=Y_{t}^{t,x},\left( t,x\right) \in \lbrack 0,T]\times {\bf
R}^{d}
\]

is the unique solution of the following SPDE: $0\leq t\leq T,$
%The following SPDE can be associated with BDSDE
\[
  u(t,x)=u(T,x)+\int_t^T [\L u(s,x)+f(x,u(s,x),(\sigma \nabla u)(s,x))]\d s -\int_t^T g(x,u(s,x),(\sigma \nabla u)(s,x))\d B_s.
\]
Note that $u(t,x)$ depends on $B(\cdot)$ indeed.

\begin{example}$f \equiv 0$, $g \equiv 1$,
\[
  u(t,x)=u(T,x)+\int_t^T \L u(s,x) \d s+ \int_t^T \d B_s.
\]
$W_t$ itself is a forward stochastic differential equation, the SDE is
\[
  X^{t,x}_s = x + \int_t^s \d W_t, \ \ \ \ t\leq s.
\]
\end{example}

%%%%%%%%%%%%%%%%%%%%%%%%%%%%%%%%%%

\subsection{Example and Simulation}

The structure of solution is interesting. Note that the collection
$\{\mathcal{F}_{t},\,t\in \lbrack 0,T]\}$ is neither increasing nor
decreasing, and it does not constitute filtration.

\begin{ex}
\begin{eqnarray*}
-dY_{t} &=&Z_{t}dB_{t}-Z_{t}dW_{t}, \\
Y_{T} &=&W_{T}.
\end{eqnarray*}

The solution is%
\[
Y_{t}=(B_{T}-B_{t})+W_{t}.
\]

We usually apply binomial tree model to simulate Brownian motion.
$W_t$ is a forward binomial tree and $B_T-B_t$ is a backward
binomial tree. Then the coupled binomial tree is a tetrahedron. It
is could be illustrated that all the paths $(t, W_t, B_T-B_t)$ are
within a tetrahedron.

%\begin{lstlisting}
%[Wt, Bt, t]=tetrahedron(1001, 10);
%plot(t,Wt);xlabel('t');ylabel('Wt');title('(t,Wt)')
%plot(t,repmat(Bt(N,:),N,1)-Bt);xlabel('t');ylabel('BT-Bt');title('(t,BT-Bt)')
%plot3(t, Wt, repmat(Bt(N,:),N,1)-Bt);
%\end{lstlisting}

$(W_t, B_T-B_t)$ is a coupled Brownian motion.

Figure (\ref{t_Wt}) illustrates $W_t$, Figure (\ref{t_Bt})
illustrates $B_T-B_t$, and Figure (\ref{t_Wt_Bt}) illustrates $(W_t,
B_T-B_t)$. The tetrahedron is big and the paths are concentrated by
central limit theorem.

\begin{figure}[!htbp]
  \caption{$(t,W_t)$\label{t_Wt}}
  \img{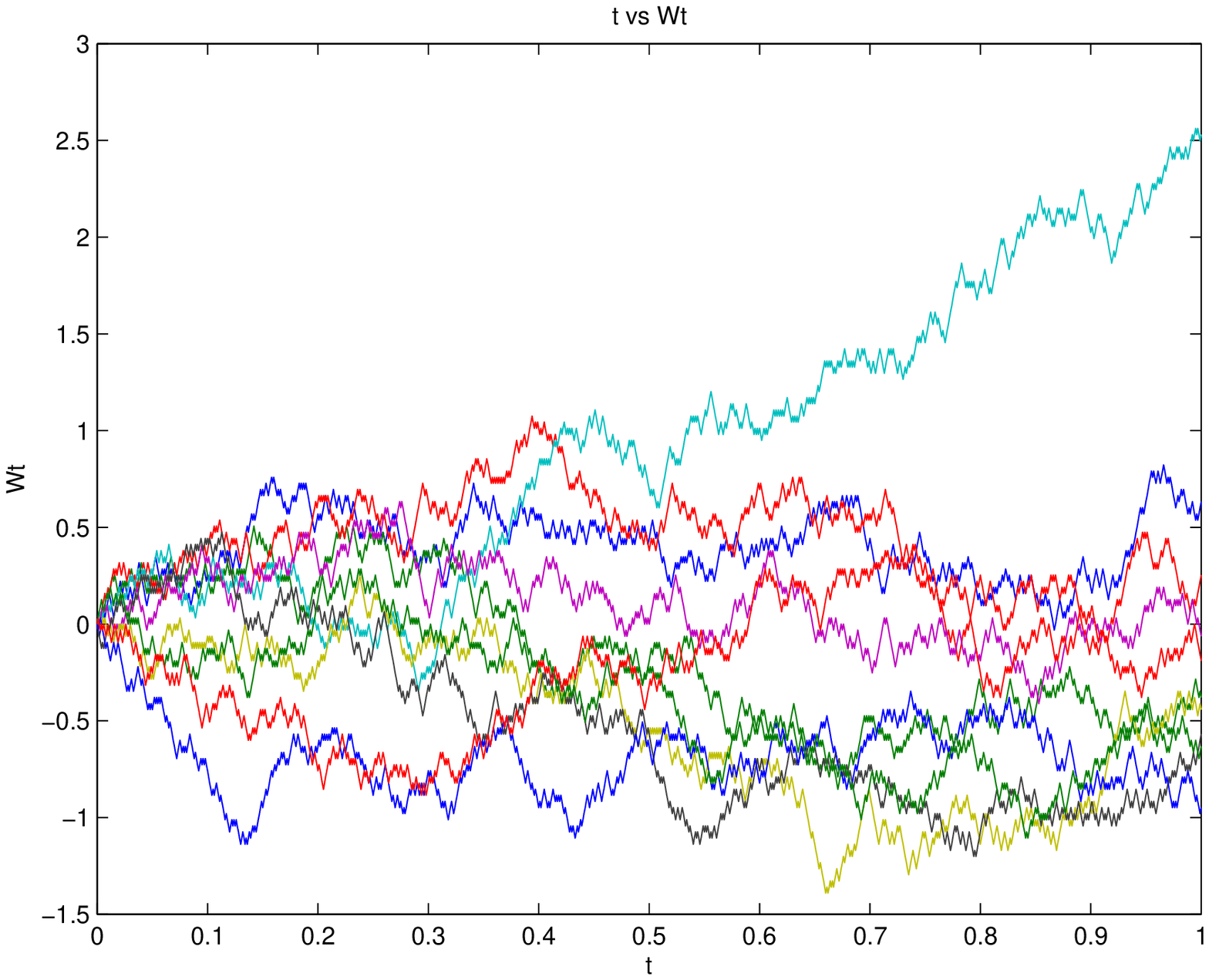}

  \caption{$(t,B_T-B_t)$\label{t_Bt}}
  \img{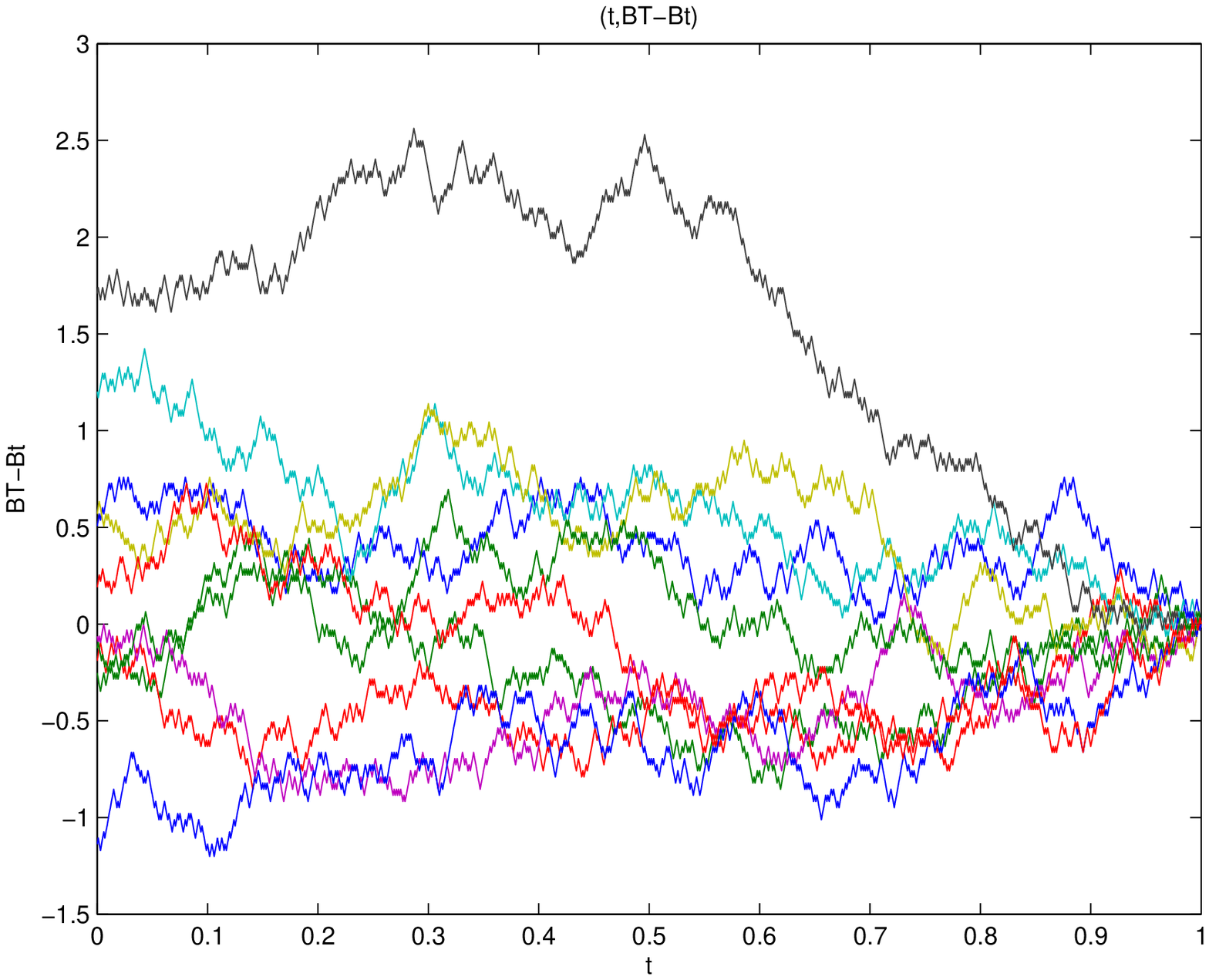}

  \caption{$(t,W_t,B_T-B_t)$\label{t_Wt_Bt}}
  \img{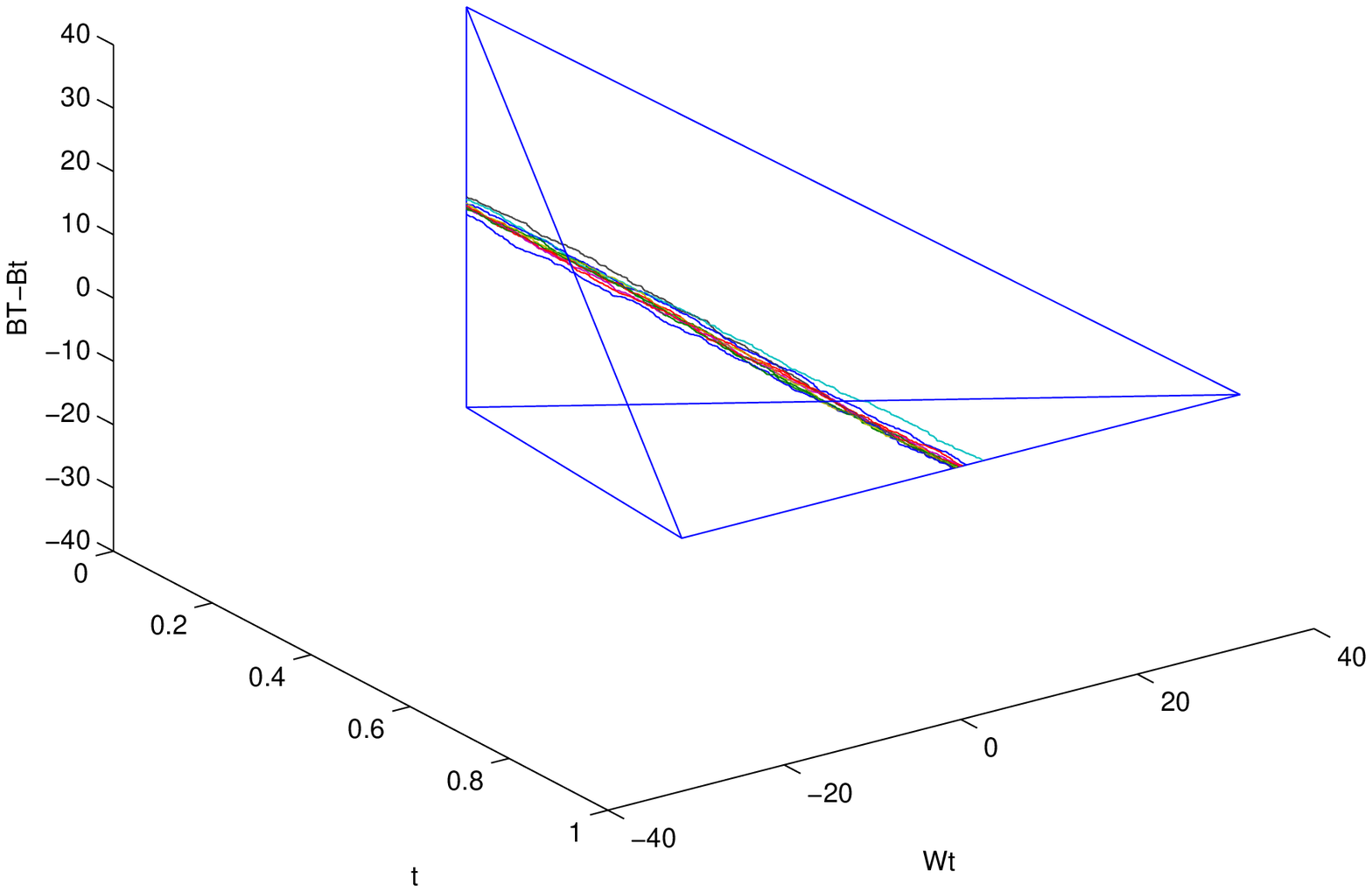}
\end{figure}

\end{ex}

\section{Main Results: Convergence Results for Discrete BDSDEs}
\label{chap:chap4}

%%%%%%%%%%%%%%%%%%%%%%%%%%%%%%%%%%%%%%%%%%%%%%%
%%%%%%%%%%%%%%%%%%%%%%%%%%%%%%%%%%%%%%%%%%%%%%%
%%%%%%%%%%%%%%%%%%%%%%%%%%%%%%%%%%%%%%%%%%%%%%%

\subsection{Convergence of The Solution for Discrete BDSDEs }
We consider the discrete terminal condition is $y_{n}^{n}:=\xi
^{n}=\Phi ((W_{j}^{n}),0\leq j\leq n)$, which is
$\mathcal{G}_{nn}^{n}$-measurable random variable, for the discrete
case. Firstly, for the scheme $(\ref{equation:discreteBDSDE})$ of
BDSDE, if we construct the processes: $$y_{t}^{n}=y_{[t/\delta
]}^{n},z_{t}^{n}=z_{[t/\delta ]}^{n}, 0\leq t\leq T$$ then the
convergence between $(y_{t}^{n},z_{t}^{n})$ to $(y_{t},z_{t})$ can
be derived in the same way as Donsker-Type theorem for BSDEs , by
(P.Briand, B. Delyon and J. M\'{e}min. (2001)\cite{PBJ}),

\textbf{Assumption} $(H.2)$ $\xi $ is $\mathcal{F}_{t}$-measurable
and, for all $n,\xi ^{n}$ is $\mathcal{G}_{nn}^{n}$-measurable s.t.
$$\mathbb{E}[\xi ^{2}] +\sup \limits_{n} \mathbb{E} [(\xi ^{n})^{2}]
<\infty$$

\textbf{Assumption} $(H.3)$ $\xi ^{n}$ converges to $\xi $ in
$L^{1}$ as $ n\rightarrow \infty.$

\bth \label{thm:convergence1}
If the assumptions $(H.1)$, $(H.2)$
and $(H.3)$ hold. Let us consider the scaled random walks
$B_{n},W_{n}$, if $B_{n}\rightarrow B,W_{n}\rightarrow W$ as
$n\rightarrow\infty$ in the sense of that $$\lim
\limits_{n\rightarrow \infty }\sup \limits_{0\leq t\leq T}
\left\vert B_{t}-B_{t}^{n}\right\vert =0\qquad \textrm{in }
\mathbf{P},$$ and $$\lim \limits_{n\rightarrow \infty }\sup
\limits_{0\leq t\leq T} \left\vert W_{t}-W_{t}^{n}\right\vert
=0\qquad \textrm{in } \mathbf{P},$$ then we have
$(y^{n},z^{n})\rightarrow (y,z)$ in the following sense:

\begin{equation} \label{equation:convergence1}
\lim \limits_{n\rightarrow \infty }\left\{ \sup \limits_{0\leq t\leq
T} \left\vert y_{t}^{n}-y_{t}\right\vert
^{2}+\int\limits_{0}^{T}\left\vert z_{s}^{n}-z_{s}\right\vert
^{2}ds\right\} =0\qquad \textrm{in } \mathbf{P}.
\end{equation}
\eth

\textbf{Method for the proof.} The key point is to use the following
decomposition

\begin{equation} \label{equation:decompositonY}
Y^{n}-Y=(Y^{n}-Y^{n,p})+(Y^{n,p}-Y^{\infty ,p})+(Y^{\infty ,p}-Y),
\end{equation}

\begin{equation} \label{equation:decompositonZ}
Z^{n}-Z=(Z^{n}-Z^{n,p})+(Z^{n,p}-Z^{\infty ,p})+(Z^{\infty ,p}-Z),
\end{equation}
where the superscript $p$ stands for the approximation of the
solution to the BDSDE via the Picard method. More precisely, we set
$Y^{\infty ,0}=0,Z^{\infty ,0}=0,Y^{n,0}=0,Z^{n,0}=0$ and define
$(Y^{\infty ,p+1},Z^{\infty ,p+1})$ as the solution of the BDSDE
\begin{equation} \label{equation:continuousPicard}
Y_{t}^{\infty ,p+1}=\xi +\int\limits_{t}^{T}f(Y_{s}^{\infty
,p},Z_{s}^{\infty ,p})ds+\int\limits_{t}^{T}g(Y_{s}^{\infty
,p},Z_{s}^{\infty
,p})dB_{s}-\int\limits_{t}^{T}Z_{s}^{\infty,p+1}dW_{s}, \quad0\leq
t\leq T.
\end{equation}

($(Y^{\infty ,p+1},Z^{\infty ,p+1})$ is solution of a BDSDE with
non-dependent but random coefficients) and similarly
\begin{eqnarray} \label{equation:discretePicard}
y_{k}^{n,p+1}&=&y_{k+1}^{n,p+1}+f(y_{k}^{n,p},z_{k}^{n,p})\delta
+g(y_{k+1}^{n,p},z_{k+1}^{n,p})\sqrt{\delta }\varepsilon
_{k+1}-z_{k}^{n,p+1} \sqrt{\delta }\beta _{k+1},\ k=n-1,...,0,{}\nonumber\\
y_n^{n,p+1}&=&\xi ^{n}
\end{eqnarray}

In order to define the discrete processes on $[0,T]$ we set for
$0\leq t\leq T$, $Y_{t}^{n,p}=y_{[t/\delta ]}^{n,p}$ and
$Z_{t}^{n,p}=z_{[t/\delta ]}^{n,p}$ so that $Y^{n,p}$ is
c\`{a}dl\`{a}g and $Z^{n,p}$ c\`{a}gl\`{a}d (c\`{a}dl\`{a}g means
right continuous with left limits and c\`{a}gl\`{a}d means left
continuous with right limits).

We shall prove in Lemma \ref{lem:convergence of discrete} that the
convergence of $(Y^{n,p},Z^{n,p})$ to $ (Y^{n},Z^{n})$ is uniform in
$n$ for the classical norm used for BDSDEs which is stronger than
the convergence in the sense of $(\ref{equation:convergence1})$;
this part is standard manipulations.

We shall prove that for any $p$, the convergence of
$(Y^{n,p},Z^{n,p})$ to $ (Y^{\infty ,p},Z^{\infty ,p})$ holds in the
sense of $(\ref{equation:convergence1})$; this is the difficult part
of the proof, and we shall need the results of section
\ref{subsubsection :Convergence of filtrations}.

\subsubsection {Convergence of Filtrations}
\label{subsubsection :Convergence of filtrations}
Let us consider a
sequence of c\`{a}dl\`{a}g processes $W^{n}=(W_{t}^{n})_{0 \leq
t\leq T}$ and $W=(W_{t})_{0\leq t\leq T}$ a Brownian motion, all
defined on the same probability space $(\Omega
,\mathcal{F},\mathbf{P})$; $T$ is finite. We denote by
$(\mathcal{G}_{t}^{n})$  (resp. $(\mathcal{G}_{t})$) the right
continuous filtration s.t.
$\sigma(W^{n})\subset\mathcal{G}^{n}_{t}$(resp.$\sigma(W)\subset\mathcal{G}_{t}$).
Let us consider finally a sequence $X^{n}$ of
$\mathcal{G}_{T}^{n}$-measurable integrable random variables, and
$X$ an $\mathcal{G}_{T}$-measurable integrable random variable,
together with the c\`{a}dl\`{a}g martingales
$$M_{t}^{n}=\mathbb{E}[X^{n}\mid
\mathcal{G}_{t}^{n}],\qquad M_{t}=\mathbb{E}[X\mid
\mathcal{G}_{t}]$$

We denote by $[M^{n},M^{n}]$ (resp. $[M,M]$) the quadratic variation
of $ M^{n}$ (resp. $M$) and by $[M^{n},W^{n}]$ (resp. $[M,W]$) the
cross variation of $M^{n}$ and $W^{n}$ (resp. $M$ and $W$).

\bth \label{thm :Convergence of filtrations}
Let us consider the
following assumptions

(A1) for each $n$, $W^{n}$ is a square integrable
$\mathcal{G}^{n}$-martingale with independent increments;

(A2) $W^{n}\rightarrow W$ in probability for the topology of uniform
convergence of c\`{a}dl\`{a}g processes indexed by $t\in [0,T]$;

(A3) a. $\mathbb{E}\left\vert X^{2}\right\vert +\sup \limits_{n}
\mathbb{E} \left\vert (X^{n})^{2}\right\vert <\infty $

\ \ \ \ b. $\mathbb{E}[\left\vert X^{n}-X\right\vert ]\rightarrow
0;$

Then, if conditions (A1) to (A3) are satisfied, we get
$$(W^{n},M^{n},[M^{n},M^{n}],[M^{n},W^{n}])\rightarrow
(W,M,[M,M],[M,W]) \qquad \textrm{in } \mathbf{P}$$ for the topology
of uniform convergence on $[0,T]$. Moreover, for each $t\in [0,T]$,
for each $0<\delta <1$,
$$(W_{t}^{n},M_{t}^{n},[M^{n},M^{n}]_{t}^{1/2},[M^{n},W^{n}]_{t}^{1/2})%
\rightarrow (W_{t},M_{t},[M,M]_{t}^{1/2},[M,W]_{t}^{1/2})\textrm {
in } L^{1+\delta }(\Omega ,\mathcal{F},\mathbf{P}).$$

\eth

\begin{corollary} Let $W$ and $W^{n}$, $n\in \mathbb{N}^{\ast }$, be
the standard Brownian motion and the random walks of Theorem
\ref{thm:convergence1} Let us consider, on the same space, $X$ and
$X^{n}$ satisfying the assumption (A3) of Theorem \ref{thm
:Convergence of filtrations}.

Then there exists a sequence $(Z_{t}^{n})_{0\leq t\leq T}$ of
$\mathcal{G} _{\cdot }^{n}$-progressively measurable processes, and
an $\mathcal{G}_{\cdot }$-progressively measurable process
$(Z_{t})_{0\leq t\leq T}$ such that: for all $t\in
[0,T]$,$$M_{t}^{n}=\mathbb{E}[X^{n}]+\int\limits_{0}^{t}Z_{s}^{n}dW_{s}^{n},\qquad
M_{t}=\mathbb{E}[X]+ \int\limits_{0}^{t}Z_{s}dW_{s}$$ and
$$\int\limits_{0}^{T}(Z_{t}^{n}-Z_{t})^{2}ds\rightarrow 0
\qquad \textrm{in } \mathbf{P}.$$

Moreover, if $0<\delta <1$, $Z^{n}$ converges to $Z$ in the space $
L^{1+\delta }(\Omega \times [0,T],\mathcal{F}\times
\mathcal{B}([0,T]),\mathbf{P}\mathscr{\otimes \lambda })$ where
$\lambda $\ denotes the Lebesgue measure on
$([0,T],\mathcal{B([}0,T\mathcal{]))}$.
\end{corollary}

\subsubsection{Proof of Theorem \ref{thm:convergence1}}

Equations
$(\ref{equation:decompositonY}),(\ref{equation:decompositonZ})$ with
the following lemma proved in appendix A.

\blm \label{lem:convergence of discrete}
Here we need to assume that
$$\lim \limits_{n\rightarrow \infty }\sup \limits_{0\leq t\leq T}
\left\vert B_{t}-B_{t}^{n}\right\vert =0\qquad \textrm{in }
\mathbf{P}.$$ With the notations following
$(\ref{equation:continuousPicard}),(\ref{equation:discretePicard}),$
$$\sup \limits_{n}\mathbb{E}[\sup \limits_{0\leq t\leq T}\left\vert
Y_{t}^{n}-Y_{t}^{n,p}\right\vert ^{2}+\int\limits_{0}^{T}\left\vert
Z_{s}^{n}-Z_{s}^{n,p}\right\vert ^{2}ds ]\rightarrow 0,\ as\
p\rightarrow \infty$$ imply that it remains to prove the convergence
to zero of the process $ Y^{n,p}-Y^{\infty ,p}$and
$Z^{n,p}-Z^{\infty ,p}.$ This will be done by induction on $p.$ For
sake of clarity, we drop the superscript $p,$set the time in
subscript and write everything in continuous time, so that $
(\ref{equation:continuousPicard}),(\ref{equation:discretePicard})$
become
$$Y_{t}^{\prime }=\xi +\int\limits_{t}^{T}f(Y_{s},Z_{s})ds+\int
\limits_{t}^{T}g(Y_{s},Z_{s})dB_{s}-\int\limits_{t}^{T}Z_{s}^{\prime
}dW_{s},\qquad0\leq t\leq T.$$

$$Y_{t}^{\prime n}=\xi
^{n}+\int\limits_{t}^{T}f(Y_{s-}^{n},Z_{s}^{n})dA_{s}^{n}+\int
\limits_{t}^{T}g(Y_{s-}^{n},Z_{s}^{n})dB_{s}^{n}-\int\limits_{t}^{T}Z_{s}^{
\prime n}dW_{s}^{n},\qquad0\leq t\leq T.$$ where
$A_{s}^{n}=[s/\delta]\delta $ and $Y_{-}$ denotes the c\`{a}gl\`{a}d
process associated to $Y$. The assumption is that
$\{Y_{t}^{n},Z_{t}^{n} \}_{0\leq t\leq T}$ converges to
$\{Y_{t},Z_{t}\}_{0\leq t\leq T}$ in sense of
$(\ref{equation:convergence1})$ and we have to prove that
$\{Y_{t}^{\prime n},Z_{t}^{\prime n}\}_{0\leq t\leq T}$ converges to
$\{Y_{t}^{\prime },Z_{t}^{\prime }\}_{0\leq t\leq T}$ in the same
sense. \elm

According to the Peng and Pardoux's paper \cite{PP}, we define the
filtration $( \mathcal{G}_{t})_{0\leq t\leq T}$ by
$$\mathcal{G}_{t}\doteq\mathcal{F}_{t}^{W}\vee \mathcal{F}_{T}^{B}$$
and the $\mathcal{G}_{t}$-square integrable martingale
$$M_{t}=\mathbb{E}^{\mathcal{G}_{t}}[\xi
+\int\limits_{0}^{T}f(Y_{s},Z_{s})ds+\int
\limits_{0}^{T}g(Y_{s},Z_{s})dB_{s}],\qquad0\leq t\leq T.$$

Then there exists $\mathcal{G}_{t}$-progressively measurable process
$ \{Z_{t}^{\prime }\}$ such that
$$\mathbb{E}\int\limits_{0}^{T}\left\vert Z_{t}^{\prime }\right\vert
^{2}dt<\infty$$
$$M_{t}=M_{0}+\int\limits_{0}^{t}Z_{s}^{\prime }dW_{s},\qquad0\leq t\leq
T.$$

On the other hand, the process, defined by
\begin{equation}
M_{t}^{n}=Y_{t}^{\prime
n}+\int\limits_{0}^{t}f(Y_{s-}^{n},Z_{s}^{n})dA_{s}^{n}+\int
\limits_{0}^{t}g(Y_{s-}^{n},Z_{s}^{n})dB_{s}^{n},\qquad0\leq t\leq
T,
\end{equation}
satisfies
\begin{equation} \label{representation of Martingale}
M_{t}^{n}=M_{0}^{n}+\int\limits_{0}^{t}Z_{s}^{\prime n}dW_{s}^{n}.
\end{equation}

Hence $M^{n}$ is an $\mathcal{F}_{\cdot }^{n}$-martingale and, since
$ Y_{T}^{n}=\xi ^{n},$
\begin{equation} \label{Martingale}
M_{t}^{n}=\mathbb{E}[M_{T}^{n}\mid \mathcal{G}_{t}^{n}],\qquad
M_{T}^{n}=Y_{t}^{n}+\int
\limits_{0}^{T}f(Y_{s-}^{n},Z_{s}^{n})dA_{s}^{n}+\int
\limits_{0}^{T}g(Y_{s-}^{n},Z_{s}^{n})dB_{s}^{n}.
\end{equation}

If we want to apply Corollary, we have to prove the $L^{1}$
convergence of $M_{T}^{n}$. But since $Y^{n}$ and $Z^{n}$ are
piecewise constant, we have
$$
\begin{array}{llllll}
&\hspace{2mm}&\left\vert
M_{T}^{n}-Y_{T}-\int\limits_{0}^{T}f(Y_{s},Z_{s})ds-\int
\limits_{0}^{T}g(Y_{s},Z_{s})dB_{s}\right\vert\\
&\leq &\left\vert Y_{T}^{n}-Y_{T}\right\vert
+\int\limits_{0}^{T}\left\vert
f(Y_{s}^{n},Z_{s}^{n})-f(Y_{s},Z_{s})\right\vert ds+\left\vert
\int\limits_{0}^{T}g(Y_{s}^{n},Z_{s}^{n})dB_{s}^{n}-\int
\limits_{0}^{T}g(Y_{s},Z_{s})dB_{s}\right\vert \\
&\leq &(1+KT)\sup \limits_{0\leq t\leq T}\left\vert
Y_{t}^{n}-Y_{t}\right\vert +K\int\limits_{0}^{T}\left\vert
Z_{s}^{n}-Z_{s}\right\vert ds+ \left\vert
\int\limits_{0}^{T}(g(Y_{s}^{n},Z_{s}^{n})-g(Y_{s},Z_{s}))dB_{s}\right\vert\\
&\hspace{2mm}&+\left\vert
\int\limits_{0}^{T}g(Y_{s}^{n},Z_{s}^{n})d(B_{s}^{n}-B_{s})\right\vert  \\
&\leq &(1+KT)\sup \limits_{0\leq t\leq T}\left\vert
Y_{t}^{n}-Y_{t}\right\vert +K\int\limits_{0}^{T}\left\vert
Z_{s}^{n}-Z_{s}\right\vert ds+K\sqrt{T}\sup \limits_{0\leq t\leq T}
\left\vert Y_{t}^{n}-Y_{t}\right\vert \\
&\hspace{2mm}&+\left\vert \int\limits_{0}^{T}(K\left\vert
Y_{s}^{n}\right\vert +\alpha \left\vert Z_{s}^{n}\right\vert
+\left\vert g(0,0)\right\vert )d(B_{s}^{n}-B_{s})\right\vert
\end{array}
$$
which tends to zero in probability and then in $L^{1}$ by
$L^{2}$-bounded. This and equations $(\ref{representation of
Martingale}),$ $(\ref{Martingale})$, imply together with Corollary
that $M^{n}$ converges to
$$M_{t}=\mathbb{E}^{\mathcal{G}_{t}}[\xi
+\int\limits_{0}^{T}f(Y_{s},Z_{s})ds+\int
\limits_{0}^{T}g(Y_{s},Z_{s})dB_{s}]=Y_{t}^{\prime
}+\int\limits_{0}^{t}f(Y_{s},Z_{s})ds+\int
\limits_{0}^{t}g(Y_{s},Z_{s})dB_{s}$$ in the sense that $$\sup
\limits_{0\leq t\leq T}\left\vert M_{t}^{n}-M_{t}\right\vert
+\int\limits_{0}^{T}\left\vert Z_{s}^{\prime n}-Z_{s}\right\vert
^{2}ds\rightarrow 0\qquad \textrm{in } \mathbf{P}.$$

Since we want to prove that $$ \sup \limits_{0\leq t\leq
T}\left\vert Y_{t}^{\prime n}-Y_{t}^{\prime }\right\vert
+\int\limits_{0}^{T}\left\vert Z_{s}^{\prime n}-Z_{s}\right\vert
^{2}ds\rightarrow 0\qquad \textrm{in } \mathbf{P},$$ it remain only
to demonstrate
$$\sup \limits_{0\leq t\leq T}\left\vert
\int\limits_{0}^{t}f(Y_{s}^{n},Z_{s}^{n})dA_{s}^{n}-\int
\limits_{0}^{t}f(Y_{s},Z_{s})ds+\int
\limits_{0}^{t}g(Y_{s-}^{n},Z_{s}^{n})dB_{s}^{n}-\int
\limits_{0}^{t}g(Y_{s},Z_{s})dB_{s}\right\vert \rightarrow 0\quad
\textrm{in } \mathbf{P}.$$

This is true since we have just proved the convergence of $
\int\limits_{0}^{T}\left\vert
f(Y_{s}^{n},Z_{s}^{n})-f(Y_{s},Z_{s})\right\vert ds$ to zero in
probability and since the jumps of $t\rightarrow
\int\limits_{0}^{t}f(Y_{s}^{n},Z_{s}^{n})dA_{s}^{n}$ tends to zero
according to $$\sup \limits_{0\leq t\leq T}\left\vert
\int\limits_{0}^{t}Z_{s}^{n}dA_{s}^{n}-\int\limits_{0}^{t}Z_{s}ds\right
\vert \rightarrow 0\quad \textrm{in } \mathbf{P},\quad\sup
\limits_{0\leq t\leq T}\left\vert
\int\limits_{0}^{t}(Z_{s}^{n})^{2}dA_{s}^{n}-\int
\limits_{0}^{t}Z_{s}^{2}ds\right\vert \rightarrow 0\quad \textrm{in
} \mathbf{P},$$ while we also have proved the convergence of
$\left\vert \int\limits_{0}^{T} g(Y_{s}^{n},Z_{s}^{n})dB_{s}^{n}-
\int\limits_{0}^{T} g(Y_{s},Z_{s})dB_{s}\right\vert$ to zero in
probability.

\subsection{Convergence of Modified Solution}
\bth \label{thm:convergence2} If the assumptions $(H.1)$, $(H.2)$
and $(H.3)$ hold. We also assume that $$\lim \limits_{n\rightarrow
\infty }\sup \limits_{0\leq t\leq T} \left\vert
B_{t}-B_{t}^{n}\right\vert =0\qquad \textrm{in } \mathbf{P}.$$ Then
the discrete solutions $\{(\overline{y}^{n},\overline{z}
^{n})\}_{n=1}^{\infty }$ under the scheme
$(\ref{equation:modifieddiscreteBDSDE})$ converge to the solution
$(y,z)$ of $(\ref{equation:BDSDE})$ in the following senses:
\begin{equation} \label{equation:convergence2}
\lim \limits_{n\rightarrow \infty }\left\{\sup \limits_{0\leq t\leq
T} \left\vert \overline{y}_{t}^{n}-y_{t}\right\vert
^{2}+\int\limits_{0}^{T}\left\vert
\overline{z}_{s}^{n}-z_{s}\right\vert ^{2}ds\right\} =0\qquad
\textrm{in } \mathbf{P}.
\end{equation}
\eth

This can be derived from the convergence
$(\ref{equation:convergence1}).$

For the convergence of this scheme, we must consider the following
estimates under the following:

\textbf{Assumption} $\xi ^{n}\in L^{2}(\xi _{n}^{n})$,
$\mathbb{E}[\sum\limits_{j=0}^{n}\left\vert
f(0,0)\right\vert^{2}]<\infty$,
$\mathbb{E}[\sum\limits_{j=0}^{n}\left\vert
g(0,0)\right\vert^{2}]<\infty .$

For this reason, we need the following Gronwall type lemma, which is
proved in \cite{MPX}.

\blm \label{lem:Gronwall}
Let us consider $a,b,\alpha $ positive
constant, $b\delta <1$ and a sequence $(v_{j})_{j=1,...,n}$ of
positive numbers such that, for every $k$
\begin{equation}
v_{j}+\alpha \leq a+b\delta \sum\limits_{i=1}^{j} v_{i},
\end{equation}
then
\begin{equation}
\sup \limits_{j\leq n}v_{j}+\alpha \leq a\varepsilon_\delta(b),
\end{equation}
where $\varepsilon _{\delta }(b)$ is the convergent sequence:
\begin{equation}
\varepsilon _{\delta }(b)=1+\sum\limits_{p=1}^{\infty
}\frac{b^{p}}{p!} (1+\delta )\cdot \cdot \cdot (1+(p-1)\delta )
\end{equation}
which is decreasing in $\delta $ and tends to $e^{b}$ as $\delta
\rightarrow \infty.$
\elm

\blm \label{lem:control of modified solution}
 We assume that $\delta $ is small enough such that
$(1+2K+7K^2)\delta<1$. Then
\begin{equation} \label{equation:control of modified solution}
\sup \limits_{j}\mathbb{E}\left\vert \overline{y}_{j}^{n}\right\vert
^{2}+\delta\mathbb{E}[\sum\limits_{j=0}^{n}\left\vert
\overline{z}_{j}^{n}\right\vert ^{2}]\leq C_{\xi
^{n},f,g}e^{(1+2K+7K^2)T}
\end{equation}
where $C_{\xi ^{n},f,g}=\left\vert f(0,0)\right\vert^2+3\left\vert
g(0,0)\right\vert^2+(1+K\delta +3K^{2}\delta+3\alpha^2\sqrt{\delta}
)\mathbb{E}\left\vert \xi^n\right\vert^2.$
\elm

\textbf{Proof.} By explicit scheme, we have
$$
\overline{y}_{j}^{n}=\overline{y}_{j+1}^{n}+f(\mathbb{E}[\overline{y}_{j+1}^{n}\mid
\mathcal{G}_{jj}^{n}],\overline{z}_{j}^{n})\delta
+g(\overline{y}_{j+1}^{n}, \overline{z}_{j+1}^{n})\sqrt{\delta
}\varepsilon _{j+1}-\overline{z}_{j}^{n} \sqrt{\delta }\beta _{j+1}.
$$

We then have
\begin{equation}
\begin{array}{ll}
\mathbb{E}\left\vert \overline{y}_{j}^{n}\right\vert
^{2}-\mathbb{E}\left\vert \overline{y} _{j+1}^{n}\right\vert
^{2}=-\mathbb{E}\left\vert \overline{z}_{j}^{n}\right\vert
^{2}\delta
+\mathbb{E}[g(\overline{y}_{j+1}^{n},\overline{z}_{j+1}^{n})]^{2}\delta
-\mathbb{E}\left\vert f(E[\overline{y}_{j+1}^{n}\mid
\mathcal{G}_{jj}^{n}],\overline{z }_{j}^{n})\right\vert ^{2}\delta
^{2}\\
\hspace{37mm}+2\mathbb{E}[\overline{y}_{j}^{n}\cdot
f(\mathbb{E}[\overline{y} _{j+1}^{n}\mid
\mathcal{G}_{jj}^{n}],\overline{z}_{j}^{n})]\delta.
\end{array}
\end{equation}

Taking sum for $j=i,...,n-1$ yields

\begin{eqnarray*}
\mathbb{E}\left\vert \overline{y}_{i}^{n}\right\vert ^{2}&\leq&
\mathbb{E}\left\vert \xi ^{n}\right\vert
^{2}-\sum\limits_{j=i}^{n-1}\mathbb{E}\left\vert \overline{z}
_{j}^{n}\right\vert ^{2}\delta +2\delta
\mathbb{E}\sum\limits_{j=i}^{n-1}\{\left\vert
\overline{y}_{j}^{n}\right\vert (\left\vert f(0,0)\right\vert
+K\left\vert \mathbb{E}[ \overline{y}_{j+1}^{n}\mid
\mathcal{G}_{jj}^{n}]\right\vert
+K\left\vert \overline{z}_{j}^{n}\right\vert )\}\\
& &+\delta \mathbb{E}\sum\limits_{j=i}^{n-1}(\left\vert
g(0,0)\right\vert +K\left\vert \overline{y}_{j+1}^{n}\right\vert
+\alpha \left\vert \overline{z}_{j+1}^{n}\right\vert )^{2}.
\end{eqnarray*}

Since the second last term is dominated by $$\delta
\mathbb{E}\sum\limits_{j=i}^{n-1}\{\left\vert
\overline{y}_{j}^{n}\right\vert ^{2}(1+K+4K^{2})+\left\vert
f(0,0)\right\vert ^{2}+K\left\vert \mathbb{E}[\overline{y}
_{j+1}^{n}\mid \mathcal{G}_{jj}^{n}]\right\vert
^{2}+\frac{1}{4}\left\vert \overline{z}_{j}^{n}\right\vert ^{2}\}$$

$$\leq \delta \mathbb{E}\sum\limits_{j=i}^{n-1}\{\left\vert \overline{y}
_{j}^{n}\right\vert ^{2}(1+K+4K^{2})+\left\vert f(0,0)\right\vert
^{2}+\frac{ 1}{4}\left\vert \overline{z}_{j}^{n}\right\vert
^{2}\}+K\delta \mathbb{E}\left\vert \xi ^{n}\right\vert ^{2}$$ and
the last term is dominated by
$$3\delta
\mathbb{E}\sum\limits_{j=i}^{n-1}(\left\vert g(0,0)\right\vert
^{2}+K^{2}\left\vert \overline{y}_{j+1}^{n}\right\vert ^{2}+\alpha
^{2}\left\vert \overline{z}_{j+1}^{n}\right\vert
^{2})\hspace{50mm}$$

$$\leq
3\delta \mathbb{E}\sum\limits_{j=i}^{n-1}(\left\vert
g(0,0)\right\vert ^{2}+K^{2}\left\vert
\overline{y}_{j}^{n}\right\vert ^{2}+\alpha ^{2}\left\vert
\overline{z}_{j}^{n}\right\vert
^{2})+3(K^2\delta+\alpha^2\sqrt{\delta})\mathbb{E}\left\vert \xi
^{n}\right\vert ^{2},$$ we thus have
\begin{eqnarray*}
\mathbb{E}\left\vert \overline{y}_{j}^{n}\right\vert ^{2}+\delta
(\frac{3}{4} -3\alpha
^{2})\sum\limits_{j=i}^{n-1}\mathbb{E}\left\vert \overline{z}
_{j}^{n}\right\vert ^{2} &\leq & \left\vert f(0,0)\right\vert
^{2}+3\left\vert g(0,0)\right\vert ^{2}+(1+K\delta
+3K^{2}\delta\\&&+3\alpha^2\sqrt{\delta} )\mathbb{E}\left\vert \xi
^{n}\right\vert ^{2} +(1+2K+7K^{2})\delta
\sum\limits_{j=i}^{n-1}\mathbb{E}\left\vert \overline{y}
_{j}^{n}\right\vert ^{2}.
\end{eqnarray*}

Then by Lemma \ref{lem:Gronwall}, we obtain $(\ref{equation:control
of modified solution}).$

\textbf{Proof of Theorem \ref{thm:convergence2} } The convergence of
$(y^{n},z^{n})$ to $(y,z) $ is already proved above. To prove that
of $(\overline{y}^{n},\overline{z} ^{n})$, it is sufficient to prove
$$\lim \limits_{n\rightarrow \infty }\left\{ \sup\limits_{0\leq t\leq
T} \left\vert \overline{y}_{t}^{n}-y_{t}^{n}\right\vert
^{2}+\int\limits_{0}^{T}\left\vert
\overline{z}_{s}^{n}-z_{s}^{n}\right\vert ^{2}ds\right\} =0.$$

From $(\ref{equation:discreteBDSDE})$ and
$(\ref{equation:modifieddiscreteBDSDE}),$ we have

\begin{eqnarray*}
\mathbb{E}\left\vert \overline{y}_{j}^{n}-y_{j}^{n}\right\vert ^{2}
&=&\mathbb{E}\left\vert
\overline{y}_{j+1}^{n}-y_{j+1}^{n}\right\vert ^{2}-\delta
\mathbb{E}\left\vert \overline{z}_{j}^{n}-z_{j}^{n}\right\vert
^{2}+\delta \mathbb{E}\left\vert g(
\overline{y}_{j+1}^{n},\overline{z}_{j+1}^{n})
-g(y_{j+1}^{n},z_{j+1}^{n}) \right\vert ^{2}\\
&\hspace{2mm}&-\mathbb{E}\left\vert
f(y_{j}^{n},z_{j}^{n})-f(\mathbb{E}[\overline{y} _{j+1}^{n}\mid
\mathcal{G}_{jj}^{n}],\overline{z}_{j}^{n})\right\vert ^{2}\delta
^{2} \\
&\hspace{2mm}&+2\mathbb{E}[(\overline{y}_{j}^{n}-y_{j}^{n})\cdot
(f(y_{j}^{n},z_{j}^{n})-f(\mathbb{E}[ \overline{y}_{j+1}^{n}
\mid\mathcal{G}_{jj}^{n}],\overline{z} _{j}^{n}))]\delta.
\end{eqnarray*}

We then take sum over $i$ from $j$ to $n-1.$ With $\xi
^{n}-\overline{\xi } ^{n}=0,$ we have
\begin{eqnarray*}
\mathbb{E}\left\vert \overline{y}_{i}^{n}-y_{i}^{n}\right\vert ^{2}
&\leq &-\delta \mathbb{E}\sum\limits_{j=i}^{n-1}\left\vert
\overline{z}_{j}^{n}-z_{j}^{n}\right\vert ^{2}+\delta
\sum\limits_{j=i}^{n-1}\mathbb{E}\left\vert
g(\overline{y}_{j+1}^{n},
\overline{z}_{j+1}^{n})-g(y_{j+1}^{n},z_{j+1}^{n})\right\vert
^{2}\\
&\hspace{2mm}&+2\delta
\sum\limits_{j=i}^{n-1}\mathbb{E}[(\overline{y}_{j}^{n}-y_{j}^{n})\cdot
(f(y_{j}^{n},z_{j}^{n})-f(\mathbb{E}[\overline{y}_{j+1}^{n}\mid
\mathcal{G}_{jj}^{n}],
\overline{z}_{j}^{n}))] \\
&\leq &-\delta \mathbb{E}\sum\limits_{j=i}^{n-1}\left\vert
\overline{z} _{j}^{n}-z_{j}^{n}\right\vert ^{2}+2K^{2}\delta
\sum\limits_{j=i}^{n-1}\mathbb{E}\left\vert
\overline{y}_{j+1}^{n}-y_{j+1}^{n}\right \vert ^{2}+2\alpha
^{2}\delta \sum\limits_{j=i}^{n-1}\mathbb{E}\left\vert \overline{z}
_{j+1}^{n}-z_{j+1}^{n}\right\vert ^{2}\\
&\hspace{2mm}&+2K^{2}\delta
\sum\limits_{j=i}^{n-1}\mathbb{E}\left\vert
\overline{y}_{j}^{n}-y_{j}^{n}\right\vert ^{2}+\delta
/2\sum\limits_{j=i}^{n-1}\mathbb{E}\left\vert \overline{z}
_{j}^{n}-z_{j}^{n}\right\vert ^{2} \\
&&+2K\delta \mathbb{E}\sum\limits_{j=i}^{n-1}\left\vert \overline{y}
_{i}^{n}-y_{i}^{n}\right\vert \cdot \left\vert
y_{j}^{n}-\mathbb{E}[\overline{y} _{j+1}^{n}\mid
\mathcal{G}_{jj}^{n}]\right\vert .
\end{eqnarray*}

Since $y_{j}^{n}-\mathbb{E}[\overline{y}_{j+1}^{n}\mid
\mathcal{G}_{jj}^{n}]=f(\mathbb{E}[ \overline{y}_{j+1}^{n}\mid
\mathcal{G}_{jj}^{n}],\overline{z}_{j}^{n})\delta +\sqrt{\delta
}\mathbb{E}[g(\overline{y}_{j+1}^{n},\overline{z}_{j+1}^{n})\mid
\mathcal{G}_{jj}^{n}]\varepsilon _{j+1}$, the last term is dominated
by
\begin{eqnarray*}
(K^2\delta^4+\delta^3&+&2K\delta^2)\sum\limits_{j=i}^{n-1}\mathbb{E}\left\vert
\overline{y}_{j}^{n}-y_{j}^{n}\right\vert
^{2}+\sum\limits_{j=i}^{n-1}K^2\mathbb{E}\left\vert f(\mathbb{E}[
\overline{y}_{j+1}^{n}\mid
\mathcal{G}_{jj}^{n}],\overline{z}_{j}^{n})\right\vert^2\delta^3
\\&&+\sum\limits_{j=i}^{n-1}\mathbb{E}\left\vert\mathbb{E}[g(\overline{y}_{j+1}^{n},\overline{z}_{j+1}^{n})\mid
\mathcal{G}_{jj}^{n}]\right\vert^2\delta.
\end{eqnarray*}

But with $(\ref{equation:control of modified solution})$, the second
term is bounded by $C\delta^2$, and the last term is bounded by
$$4K^2\delta\sum\limits_{j=i}^{n-1}\mathbb{E}\left\vert
\overline{y}_{j}^{n}\right\vert
^{2}+4\alpha^2\delta\sum\limits_{j=i}^{n-1}\mathbb{E}\left\vert
\overline{z}_{j}^{n}\right\vert
^{2}+4(K^2\delta+\alpha^2\sqrt{\delta})\mathbb{E}\left\vert\xi^n\right\vert^2+2\left\vert
g(0,0)\right\vert^2.$$ We thus have
\begin{eqnarray*}
\mathbb{E}\left\vert \overline{y}_{i}^{n}-y_{i}^{n}\right\vert
^{2}+(\frac{1}{2}-2\alpha^2)\delta
\mathbb{E}\sum\limits_{j=i}^{n-1}\left\vert
\overline{z}_{j}^{n}-z_{j}^{n}\right\vert ^{2} &\leq &
(4K^2\delta+K^2\delta^4+\delta^3+2K\delta^2)\sum\limits_{j=i}^{n-1}\mathbb{E}\left\vert
\overline{y}_{j}^{n}-y_{j}^{n}\right\vert ^{2}\\
&\hspace{1mm}&
+4K^2\delta\sum\limits_{j=i}^{n-1}\mathbb{E}\left\vert
\overline{y}_{j}^{n}\right\vert
^{2}+4\alpha^2\delta\sum\limits_{j=i}^{n-1}\mathbb{E}\left\vert
\overline{z}_{j}^{n}\right\vert
^{2}\\
&\hspace{1mm}&+4(K^2\delta+\alpha^2\sqrt{\delta})\mathbb{E}\left\vert\xi^n\right\vert^2+2\left\vert
g(0,0)\right\vert^2+C\delta^2.
\end{eqnarray*}
According to Lemma \ref{lem:control of modified solution} and here
providing that $g(0,0)=0$, we further have
\begin{eqnarray*}
\mathbb{E}\left\vert \overline{y}_{i}^{n}-y_{i}^{n}\right\vert
^{2}+(\frac{1}{2}-2\alpha^2)\delta
\mathbb{E}\sum\limits_{j=i}^{n-1}\left\vert
\overline{z}_{j}^{n}-z_{j}^{n}\right\vert ^{2} &\leq &
(1+2k+5k^2)\delta\sum\limits_{j=i}^{n-1}\mathbb{E}\left\vert
\overline{y}_{j}^{n}-y_{j}^{n}\right\vert ^{2}+C'\sqrt{\delta}.
\end{eqnarray*}
By Gronwall Lemma 4.2.2, we get
$$\sup\limits_{i\leq n}\mathbb{E}\left\vert \overline{y}_{i}^{n}-y_{i}^{n}\right\vert
^{2}\leq C'\sqrt{\delta} e^{(1+2k+5k^2)T}.
$$
Then these two inequalities implies (\ref{equation:convergence2}).

\section*{Appendix A.  Proof of Lemma \ref{lem:convergence of discrete}} \label{chap:appendix A.}

\renewcommand{\theequation}{A.\arabic{equation}}
  % redefine the command that creates the equation no.
\setcounter{equation}{0}
% reset counter

For the proof of this lemma we come back to the discrete notations
and we show that
\newline\textbf{Lemma A.1} \emph{There exist $\alpha>1$
and $n_0\in \mathbb{N}$ such that for all $n\geq n_0$, for all
$p\in\mathbb{N}^*$,
$$\left\Vert (y^{n,p+1}-y^{n,p},z^{n,p+1}-z^{n,p})\right\Vert_\alpha
^{2}\leq \frac{2}{3}\left\Vert
(y^{n,p}-y^{n,p-1},z^{n,p}-z^{n,p-1})\right\Vert_\alpha ^{2},$$
where, for $p\in\mathbb{N}$,
$$\left\Vert (y^{n,p+1}-y^{n,p},z^{n,p+1}-z^{n,p})\right\Vert_\alpha
^{2}:=\mathbb{E}\left [\sup\limits_{0\leq k\leq n}
\alpha^{k\delta}\left\vert y_k^{n,p+1}-y_k^{n,p}
\right\vert^2+\delta\sum\limits_{k=0}^{n-1}
\alpha^{k\delta}\left\vert z_k^{n,p+1}-z_k^{n,p} \right\vert^2\right
].$$}

\textbf{Proof.} For notational convenience, let us write $y,z$
in place of $ y^{n,p+1}-y^{n,p},$ $z^{n,p+1}-z^{n,p}$ and $u,v$ in
place of $ y^{n,p}-y^{n,p-1},$ $z^{n,p}-z^{n,p-1}$. Let us pick
$\varphi >1$ to be chosen later. With these notations in hands, we
have, for $k=0,...,n-1$, since $y_{n}=0,$
$$\varphi
^{k}y_{k}^{2}=\sum\limits_{i=k}^{n-1}\left( \varphi
^{i}y_{i}^{2}-\varphi ^{i+1}y_{i+1}^{2}\right) =(1-\varphi
)\sum\limits_{i=k}^{n-1}\varphi ^{i}y_{i}^{2}+\varphi
\sum\limits_{i=k}^{n-1}\varphi ^{i}(y_{i}^{2}-y_{i+1}^{2}).$$

We write
$y_{i}^{2}-y_{i+1}^{2}=2y_{i}(y_{i}-y_{i+1})-(y_{i}-y_{i+1})^{2},$
to use the equation $(\ref{equation:discretePicard}),$ since
\begin{eqnarray}
y_{i}-y_{i+1}&=&\delta
\{f(y_{i}^{n,p},z_{i}^{n,p})-f(y_{i}^{n,p-1},z_{i}^{n,p-1})\}+\sqrt{\delta
}
\{g(y_{i+1}^{n,p},z_{i+1}^{n,p})-g(y_{i+1}^{n,p-1},z_{i+1}^{n,p-1})\}\varepsilon
_{i+1}\nonumber
\\& &
-\sqrt{\delta }z_{i}\beta_{i+1}.
\end{eqnarray}

According to $(H.1),$ we have, for each $\nu >0,$
\begin{eqnarray*}
2y_{i}\{f(y_{i}^{n,p},z_{i}^{n,p})-f(y_{i}^{n,p-1},z_{i}^{n,p-1})\}&\leq&
2K\left\vert y_{i}\right\vert (\left\vert u_{i}\right\vert
+\left\vert v_{i}\right\vert )\\
&\leq& 2(K^{2}/\nu )y_{i}^{2}+\nu (u_{i}^{2}+v_{i}^{2}),
\end{eqnarray*}
\begin{eqnarray*}
2y_{i}\{g(y_{i+1}^{n,p},z_{i+1}^{n,p})-g(y_{i+1}^{n,p-1},z_{i+1}^{n,p-1})\}
&\leq& 2\left\vert y_{i}\right\vert (K\left\vert u_{i+1}\right\vert
+\alpha \left\vert v _{i+1}\right\vert )\\
&\leq& \lbrack (K^{2}+\alpha ^{2})/\nu]y_{i}^{2}+\nu
(u_{i+1}^{2}+v_{i+1}^{2})
\end{eqnarray*}
and moreover, $(A.1)$ implies easily that
$$\delta z_{i}^{2}\leq 3(y_{i}-y_{i+1})^{2}+6K^{2}\delta ^{2}(u
_{i}^{2}+v_{i}^{2})-6\delta (K^{2}u_{i+1}^{2}+\alpha ^{2}v
_{i+1}^{2}).$$

As a byproduct of these inequalities, we deduce that, for
$k=0,...,n-1,$

\begin{eqnarray*}
&&2\sum\limits_{i=k}^{n-1}\varphi ^{i}y_{i}(y_{i}-y_{i+1})\\
&\leq&2K^{2}(\delta /\nu)\sum\limits_{i=k}^{n-1}\varphi
^{i}y_{i}^{2}+\nu\delta \sum\limits_{i=k}^{n-1}\varphi
^{i}(u_{i}^{2}+v_{i}^{2})+[\sqrt{\delta }(K^{2}+\alpha ^{2})/\nu
]\sum\limits_{i=k}^{n-1}\varphi ^{i}y_{i}^{2}\varepsilon _{i+1} \\
&&+\nu \sqrt{\delta }\sum\limits_{i=k}^{n-1}\varphi ^{i}(u
_{i+1}^{2}+v_{i+1}^{2})\varepsilon _{i+1}-2\sqrt{\delta }
\sum\limits_{i=k}^{n-1}\varphi ^{i}y_{i}z_{i}\beta _{i+1},
\end{eqnarray*}

$$-\sum\limits_{i=k}^{n-1}\varphi ^{i}(y_{i}^{2}-y_{i+1}^{2})\leq
-(\delta /3)\sum\limits_{i=k}^{n-1}\varphi
^{i}z_{i}^{2}+2K^{2}\delta ^{2}\sum\limits_{i=k}^{n-1}\varphi
^{i}(u_{i}^{2}+v_{i}^{2})-2\delta \sum\limits_{i=k}^{n-1}\varphi
^{i}(K^{2}u_{i+1}^{2}+\alpha ^{2}v_{i+1}^{2}),$$ and, setting $ \rho
=(\nu+2K^{2}\delta )\varphi \delta ,$ we get
\begin{eqnarray*}
\varphi ^{k}y_{k}^{2}+\varphi (\delta
/3)\sum\limits_{i=k}^{n-1}\varphi ^{i}z_{i}^{2} &\leq &(1-\varphi
+2K^{2}\delta \varphi /\nu)\sum\limits_{i=k}^{n-1}\varphi
^{i}y_{i}^{2}+\rho \sum\limits_{i=k}^{n-1}\varphi ^{i}(u
_{i}^{2}+v_{i}^{2})\\
&\hspace{2mm}&+[\sqrt{\delta }\varphi (K^{2}+\alpha ^{2})/\nu
]\sum\limits_{i=k}^{n-1}\varphi ^{i}y_{i}^{2}\varepsilon _{i+1}
+\nu\varphi \sqrt{\delta }\sum\limits_{i=k}^{n-1}\varphi
^{i}(u_{i+1}^{2}+v _{i+1}^{2})\varepsilon _{i+1}\\
&\hspace{2mm}&-2\delta \varphi \sum\limits_{i=k}^{n-1}\varphi
^{i}(K^{2}u_{i+1}^{2}+\alpha ^{2}v _{i+1}^{2})-2\varphi \sqrt{\delta
}\sum\limits_{i=k}^{n-1}\varphi ^{i}y_{i}z_{i}\beta _{i+1}.
\end{eqnarray*}

Thus, if $1-\varphi +2K^{2}\delta \varphi /\nu \leq 0$, we have, for
$ k=0,...,n-1,$
\begin{equation}
\begin{array}{ll}
\varphi ^{k}y_{k}^{2}+\varphi (\delta
/3)\sum\limits_{i=k}^{n-1}\varphi ^{i}z_{i}^{2}\\
\leq\rho \sum\limits_{i=k}^{n-1}\varphi ^{i}(u _{i}^{2}+v
_{i}^{2})+[\sqrt{\delta }\varphi (K^{2}+\alpha ^{2})/\upsilon
]\sum\limits_{i=k}^{n-1}\varphi ^{i}y_{i}^{2}\varepsilon _{i+1} \\
+\nu\varphi \sqrt{\delta }\sum\limits_{i=k}^{n-1}\varphi
^{i}(u_{i+1}^{2}+v_{i+1}^{2})\varepsilon _{i+1}-2\varphi
\sqrt{\delta } \sum\limits_{i=k}^{n-1}\varphi ^{i}y_{i}z_{i}\beta
_{i+1}.
\end{array}
\end{equation}

In particular, taking the expectation of the previous inequality for
$k=0$, we get
\begin{equation}
\mathbb{E}[\sum\limits_{i=0}^{n-1}\varphi ^{i}z_{i}^{2}]\leq
3(\nu+2K^{2}\delta )\mathbb{E}[\sum\limits_{i=0}^{n-1}\varphi
^{i}(u_{i}^{2}+v_{i}^{2})].
\end{equation}

Now, coming back to $(A.2),$ we have, since $y_{n}=0,$
\begin{eqnarray*}
\sup \limits_{0\leq k\leq n}\varphi ^{k}y_{k}^{2} &\leq &\rho
\sum\limits_{i=0}^{n-1}\varphi ^{i}(u_{i}^{2}+v_{i}^{2})+2[\sqrt{
\delta }\varphi (K^{2}+\alpha ^{2})/\nu ] \sup \limits_{0\leq k\leq
n-1}\left\vert \sum\limits_{i=0}^{n-1}\varphi
^{i}y_{i}^{2}\varepsilon
_{i+1}\right\vert \\
&&+2\nu \varphi \sqrt{\delta }\sup \limits_{0\leq k\leq n-1}
\left\vert \sum\limits_{i=0}^{n-1}\varphi ^{i}(u_{i+1}^{2}+v
_{i+1}^{2})\varepsilon _{i+1}\right\vert +4\varphi \sqrt{\delta
}\sup \limits_{ 0\leq k\leq n-1}\left\vert
\sum\limits_{i=0}^{n-1}\varphi ^{i}y_{i}z_{i}\beta _{i+1}\right\vert
\end{eqnarray*}
and using Burkholder-Davis-Gundy inequality, we obtain,
\begin{eqnarray*}
\mathbb{E}[\sup \limits_{0\leq k\leq n}\varphi ^{k}y_{k}^{2}] &\leq
& \rho \mathbb{E}[\sum\limits_{i=0}^{n-1}\varphi ^{i}(u _{i}^{2}+v
_{i}^{2})]+C_{1}[ \sqrt{\delta }\varphi (K^{2}+\alpha ^{2})/\nu
]\mathbb{E}[(\sum\limits_{i=0}^{n-1}\varphi
^{2i}y_{i}^{4})^{1/2}]\\
&&+C_{2}\nu \varphi \sqrt{\delta
}\mathbb{E}[(\sum\limits_{i=0}^{n-1}\varphi ^{2i}(u_{i+1}^{2}+v
_{i+1}^{2})^{2})^{1/2}] +C_{3}\varphi \sqrt{\delta
}\mathbb{E}[(\sum\limits_{i=0}^{n-1}\varphi
^{2i}y_{i}^{2}z_{i}^{2})^{1/2}] \\
&\leq &\rho \mathbb{E}[\sum\limits_{i=0}^{n-1}\varphi ^{i}(u
_{i}^{2}+v_{i}^{2})]+C_{1}\varphi [(K^{2}+\alpha ^{2})/\nu]
\mathbb{E}[\sup \limits_{0\leq k\leq n}\varphi ^{k}y_{k}^{2}]\\
&&+C_{2}\nu \delta \mathbb{E}[\sum\limits_{i=0}^{n-1}\varphi
^{i}(u_{i+1}^{2}+v_{i+1}^{2})] +C_{3}^{2}\varphi ^{2}\delta
\mathbb{E}[\sum\limits_{i=0}^{n-1}\varphi
^{i}z_{i}^{2}]+\frac{1}{4}\mathbb{E}[\sup \limits_{0\leq k\leq
n}\varphi ^{k}y_{k}^{2}].
\end{eqnarray*}

Finally, from $(A.3)$, we get the inequality,
\begin{equation}
\mathbb{E}[\sup \limits_{0\leq k\leq n}\varphi
^{k}y_{k}^{2}+\delta\sum\limits_{i=0}^{n-1}\beta^iz_i^2]
\leq\lambda\mathbb{E}[\sup \limits_{0\leq k\leq n}\varphi
^{k}u_{k}^{2}+\delta\sum\limits_{i=0}^{n-1}\varphi^iv_i^2],
\end{equation}
where
$\lambda=\frac{\varphi(\nu+2K^2\delta)(1+3C_3^2\varphi-3C_1(K^2+\alpha^2)/\nu)+C_2\nu}{\frac{3}{4}-C_1\varphi(K^2+\alpha^2)/\nu}$
and providing that $1-\varphi+2K^2\delta\varphi/\nu\leq0$.

Firstly, we choose $\nu$ such that
$\frac{\nu(1+3C_3^2-3C_1(K^2+\alpha^2)/\nu)+C_2\nu}{\frac{3}{4}-C_1(K^2+\alpha^2)/\nu}=1/2$.
We consider only $n$ greater than $n_1$ (i.e. $K\delta<1$ and
$2K^2\delta/\nu<1$). Let us pick $\varphi$ of the form
$\gamma^\delta$ with $\gamma\geq1$. We want that
$1-\gamma^\delta+2K^2\delta\gamma^\delta/\nu\leq0$ meaning that
$\gamma\geq exp{-\delta^{-1}log(1-2K^2\delta/\nu)}$. Since
$exp\{-\delta^{-1}log(1-2K^2\delta/\nu)\}$ tends to
$exp\{2K^2\delta/\nu\}$ as $n\rightarrow\infty
(\delta\rightarrow0)$, we choose $\gamma=exp\{1+2K^2\delta/\nu\}$.
Hence, for $n$ greater than $n_2$ the condition is satisfied and
(4.18) holds for $\varphi=\gamma^\delta$. It remains to observe
that, $\nu$ and $\gamma$ being fixed as explained above, $\lambda$
converges, as $n\rightarrow\infty$, to
$\frac{\nu(1+3C_3^2-3C_1(K^2+\alpha^2)/\nu)+C_2\nu}{\frac{3}{4}-C_1(K^2+\alpha^2)/\nu}$
which is equal to $1/2$. It follows that for $n$ large enough, say
$n\geq n_0, \lambda\leq 2/3$ and
$$\mathbb{E}[\sup \limits_{0\leq k\leq n}\gamma^{k\delta}
y_{k}^{2}+\delta\sum\limits_{i=0}^{n-1}\gamma^{i\delta}z_i^2] \leq
2/3\mathbb{E}[\sup \limits_{0\leq k\leq
n}\gamma^{k\delta}u_{k}^{2}+\delta\sum\limits_{i=0}^{n-1}\gamma^{i\delta}v_i^2],$$
which concludes the proof of this technical lemma.

To complete the proof of Lemma \ref{lem:convergence of discrete}, it
remains to check that
$$\sup\limits_n \mathbb{E}\left[ \sup\limits_{0\leq k\leq n-1}|y_k^{n,1}|^2
+\delta\sum\limits_{i=0}^{n-1}|z_i^{n,1}|^2\right]$$ is finite. But
it is plain to check (using the same computations as above) that for
$n$ large enough,
$$\mathbb{E}\left[ \sup\limits_{0\leq k\leq n-1}|y_k^{n,1}|^2
+\delta\sum\limits_{i=0}^{n-1}|z_i^{n,1}|^2\right]\leq
Ce^2\left(\mathbb{E}[\xi^2]+|f(0,0)|^2+|g(0,0)|^2\right)$$ where $C$
is a universal constant.

{\bf Acknowledgements}

The authors thank Professor Shige Peng for his helpful discussion.

%\appendix
%\backmatter

%\nocite{*}
%\small

%\bibliography{all}
%\bibliographystyle{setup/jmb} % or plain
\bibliographystyle{plain}

\end{document}